\documentclass[12pt]{article}
\usepackage{amssymb}
\usepackage{amsfonts}
\usepackage{amsmath}
\usepackage[dvips]{graphicx}
\usepackage{setspace}

\begin{document}
\begin{center}
\Large Strange Integrals Derived By Elementary Complex Analysis
\end{center}

\begin{flushright}
Josh Isralowitz \\

New Jersey Institute of Technology \\

606 Cullimore Hall \\

323 Martin Luther King, Jr. Blvd.  \\

University Heights \\

Newark, NJ 07102-1982    U.S.A. \\

\verb+jbi2@njit.edu+

\end{flushright}

\begin{abstract}
In this short note, we provide an elementary complex analytic method
for converting known real integrals into numerous strange and
interesting looking real integrals.
\end{abstract}

\vspace {5mm}

 Consider the odd looking integral

\begin{equation}
 lim_{ R \to \infty}  \int_0^\pi R \hspace{1mm} e^{-R^2 cos \hspace{1mm}
2\theta}
 sin \left(R^2 sin \hspace{1mm} 2\theta - \theta \right) \, d\theta,
\end{equation}

  which typically invites responses like "what
the...!?" Despite the appearance of such an integral, however, we
provide a method to compute integrals like equation $\left(1
\right)$ using certain well known real integrals by simply modifying
the typical complex analytic method for computing real integrals.
Moreover, the method we describe provides an easy way to motivate
the beauty and simplicity of complex analysis early on in a complex
analysis course.

\vspace {3mm}

The typical complex analytic method for solving real integrals
(which is found in nearly every complex analysis book, for example
[\textbf{1}] or  [\textbf{2}]) is to consider the integral in the
complex plane where we start with a line segment on the real line,
appropriately close the line segment, integrate over the resulting
closed curve, and finally use this complex integral to obtain the
value of the original real integral. Here is an example taken from
[\textbf{1}]. If we want to compute $\int_{- \infty}^\infty
\frac{sin
 \hspace{1mm} x}{x} \, dx, $ then we can equivalently compute the imaginary
part of the principle value of $\int_{- \infty}^\infty \frac{e^{i
z}}{z} \, dx. $

\begin{figure}[h] \centering
\includegraphics{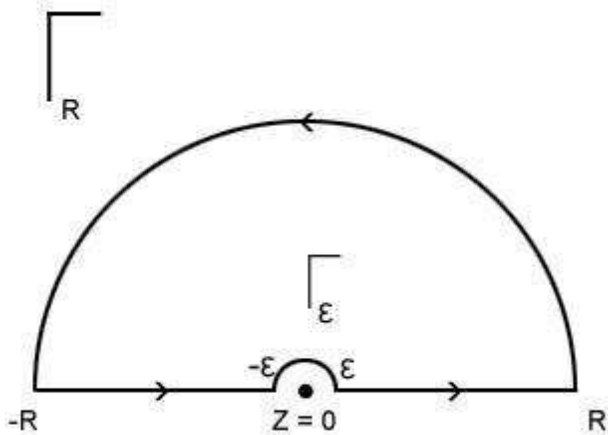}
\caption{Contour for the current example.}
\end{figure}

To do this, consider $\int_{\Gamma} \frac{e^{i z}}{z} \, dz $
 so that basic complex analytic results gives us that \[lim_{R \to \infty} \int_{\Gamma_R} \frac{e^{i
z}}{z} \, dz = 0\] and \[lim_{\epsilon \to 0} \int_{\Gamma_\epsilon}
\frac{e^{i z}}{z} \, dz = - i\pi \] where $\Gamma_R = \{z : |z| = R
$ and $Im(z) > 0 \}, \Gamma_\epsilon = \{z : |z| = \epsilon $ and
$Im(z)
> 0\} $, and $\Gamma = \Gamma_\epsilon \cup \Gamma_R \cup [-R,
- \epsilon] \cup [ \epsilon, R], $ so that
\[\int_{- \infty}^\infty \frac{sin \hspace{1mm}
  x}{x} \, dx = \pi. \]

\paragraph*{Main Examples}
Now we show how integrals like equation $\left(1 \right)$ are
derived from certain well known real integrals. We know that
$\int_{- \infty}^\infty \frac{sin \hspace{1mm}
  ax}{x} \, dx = \pi $ and that $f(z) = \frac{sin \hspace{1mm} z}{z}$ is
entire (more precisely, the function $f(z) = \frac{sin \hspace{1mm}
z}{z}$ on $z \in \mathbb{C} \backslash \{0\}$ and $f(z) = 1$ for $z
= 0$ is entire.) Therefore, we have

\begin{equation}\lim_{R \to \infty}
\int_{\Gamma_R} \frac{sin \hspace{1mm} z}{z} \, dz = - \int_{-
\infty}^\infty \frac{sin \hspace{1mm}
  x}{x} \, dx = - \pi
\end{equation}

   where again, $\Gamma_R = \{z : |z| = R $ and
$Im(z) > 0 \}.$

\vspace{5mm}

However, if we put the first integral of equation $\left(2 \right)$
in terms of real quantities, then we obtain a very odd looking
result. First notice that $z = R e^{i \theta}$ on $\Gamma_R$ so that
$dz = i R e^{i \theta} d\theta$. Therefore, directly plugging in for
$z$ and $dz$ into the first integral of equation $\left(2 \right)$
and plowing through some tedious but completely intuitive and
elementary calculations involving multiplication, the $sin$ addition
formula, the $sin/sinh$ and $cos/cosh$ identities (and only
considering the real part of the integral, since we know that the
imaginary part must be zero) gives us that
\begin{equation}
 - \lim_{R \to \infty}
\int_{\Gamma_R} \frac{sin \hspace{1mm} z}{z} \, dz = \lim_{R \to
\infty} \int_0^\pi cos \hspace{1mm} \left(R \hspace{1mm} cos
\hspace{1mm} \theta \right) sinh \hspace{1mm} \left(R \hspace{1mm}
sin \hspace{1mm} \theta \right) \, d \theta = \pi,
\end{equation}

which is indeed a very odd looking integral. Now we will show how to
use this method to derive the integral in equation $\left(1
\right).$

\vspace{5mm}

Consider the entire function $f(z) = e^{-z^2}$.  Using the same
contours and parametrization as in the previous example, we obtain
\[ - \lim_{R \to \infty} \int_{\Gamma_R} e^{-z^2} \, dz =
\int_{-\infty}^\infty e^{-x^2} \, dx = \surd \pi. \] Therefore,
substituting in for $z$ and $dz$ and again plowing through some
tedious but completely intuitive and elementary calculations
involving multiple uses of Euler's identity, multiplication, the
$sin$ addition formula (where again disregarding the imaginary part
of the integral) gives us that

 \[ - lim_{ R \to \infty}
 \int_0^\pi R \hspace{1mm} e^{-R^2 cos \hspace{1mm}  2\theta}
 sin \left(R^2 sin \hspace{1mm} 2\theta - \theta \right) \, d\theta =
 \int_{-\infty}^\infty e^{-x^2} \, dx =  \surd \pi, \] which is of
 course the integral in equation $\left(1 \right).$

\vspace{5mm} Finally, we use the same integral
$\int_{-\infty}^\infty e^{-x^2}\, dx$, but use a different curve to
close the line segment $[-R, R]$.  This time, let \[\Gamma_R = \{ t
+ i \left(t^2 - R^2 \right) : t \in [-R, R] \} \] so that again, \[
- \lim_{R \to \infty} \int_{\Gamma_R} e^{z^2} \, dz =
\int_{-\infty}^\infty e^{-x^2}\, dx = \surd \pi. \]

\begin{figure}[h] \centering
\includegraphics{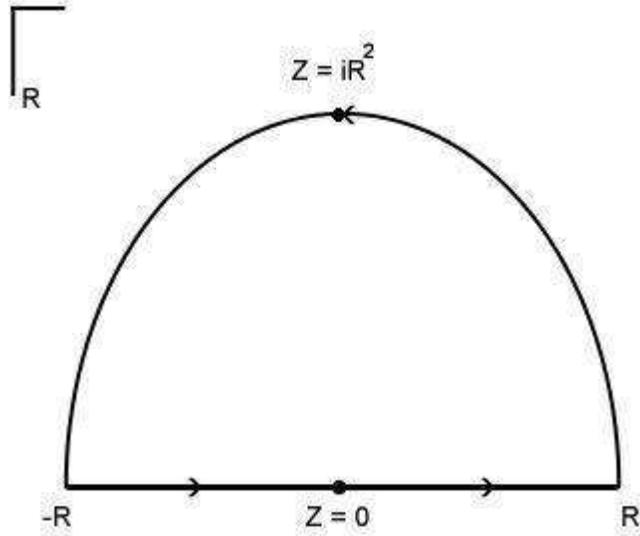}
\caption{Contour for the current example.}
\end{figure}

Thus, plugging in for $z = t + i (t^2 - R^2)$ and $dz = (1 + 2it
)dt$,  we obtain the following: \[ - \lim_{R \to \infty} \int_{-R}^R
e^{\left(t^2 - R^2 \right)^2 - t^2} \left( cos \hspace{1mm}
\left(2t^3 - tp^2 \right) + 2t sin \hspace{1mm} \left(2t^3 - tp^2
\right) \right) \, dt = \surd \pi,\] which is yet another strange
looking integral.

\paragraph*{Conclusion}
At this point, we could evaluate numerous other integrals by merely
changing the integrands or using various paths to close the line
segment $[-R, R].$ For example, we could have simply used any
polynomial with its only real roots as $\{-R, R\}$ or ellipses to
close the line segment $[-R, R]$. However, doing so involves nothing
more than reapplying our method to whatever integrand or path to
close the line segment one chooses, and so we leave this for the
interested reader to do. After all, what's more fun than deriving a
formula that is quite possibly more weird looking than
\[- \lim_{R \to \infty} \int_{-R}^R e^{(t^2 - R^2)^2
- t^2} \left( cos \hspace{1mm} \left(2t^3 - tp^2\right) + 2t sin
\hspace{1mm} \left(2t^3 - tp^2 \right) \right) \, dt ?\]

\paragraph*{Acknowledgment}
I wish to thank Prof. Jane Gilman and Prof. Jacob Sturm for their
time, effort, and useful suggestions.

\subsection*{References}

1.  Mark Ablowitz and Athanassios Fokas,  Complex Variables:
Introduction and Applications, Cambridge University Press, New York,
2003.

2.  James Brown, Reul Churchill, and Roger Verhey, Complex Variables
and Applications, McGraw-Hill, New York, 1996.

\end{document}